\documentclass[leqno,12pt]{article}
\usepackage[american,ngerman]{babel}   
\usepackage[utf8]{inputenc}

\usepackage{amsmath,amssymb,amsthm}
\usepackage{amscd}
\usepackage{setspace} 

\usepackage{enumerate}
\usepackage{array}
\usepackage{fancyhdr,ifthen}
\usepackage{tabularx}
\usepackage{multirow}

\theoremstyle{definition}
\newtheorem{definition}{Definition}
\newtheorem{remark}{Remark}

\theoremstyle{plain}

\newtheorem{lemma}{Lemma}
\newtheorem{theorem}{Theorem}
\newtheorem{corollary}{Corollary}
\newtheorem{proposition}{Proposition}

\theoremstyle{remark}
\newtheorem*{eproof}{Proof}

\newcounter{newitem}
{%
\end{list}%
}%

\newcommand{\C}{\mathbb{C}}
\newcommand{\F}{\mathbb{F}}

\newcommand{\N}{\mathbb{N}}
\newcommand{\Q}{\mathbb{Q}}

\newcommand{\Z}{\mathbb{Z}}
\newcommand{\Hcal}{\mathcal{H}}
\newcommand{\Ocal}{\mathcal{O}}
\newcommand{\Ucal}{\mathcal{U}}

\newcommand{\re}{\operatorname{Re}}

\newcommand{\trace}{\operatorname{trace}}
\newcommand{\rank}{\operatorname{rank}\,}
\newcommand{\aut}{\operatorname{Aut}}
\newcommand{\diag}{\operatorname{diag}}

\begin{document}
\thispagestyle{empty}
\vspace*{4ex}

\begin{center}
\begin{huge}
\begin{spacing}{1.0}
\textbf{On the classification of lattices over $\Q(\sqrt{-3})$, which are even unimodular $\Z$-lattices}  
\end{spacing}
\end{huge}

\vspace{1.5cm} 
by
\vspace{1.5cm}

\begin{large}
\textsc{Michael Hentschel, Aloys Krieg} 
\end{large}
\end{center}
Lehrstuhl A für Mathematik, RWTH Aachen University, D-52056 Aachen  \\
michael.hentschel@mathA.rwth-aachen.de, krieg@mathA.rwth-aachen.de 
\vspace{1cm}

\begin{center}
\begin{large}
\textsc{Gabriele Nebe}
\end{large}
\end{center} 
Lehrstuhl D für Mathematik, RWTH Aachen University, D-52056 Aachen\\
gabriele.nebe@math.rwth-aachen.de
\vspace{1.5cm}  

\newpage
\setcounter{page}{1}

\section{Introduction}
The classification of even unimodular $\Z$-lattices is explicitly known only in the cases of rank $8$, $16$ and $24$ (cf. \cite{CS}). Given an imaginary quadratic number field $K$ with discriminant $d_K$ and ring of integers $\Ocal_K$, Cohen and Resnikoff \cite{CR} showed that there exists a free $\Ocal_K$-module $\Lambda$ of $\rank r$, which is even and satisfies $\det\Lambda=(2/\sqrt{d_K})^r$ if and only if $r\equiv 0\pmod{4}$, where explicit examples can also be found in \cite{DK}. Each such $\Ocal_K$-module is an even unimodular $\Z$-lattice of $\rank 2r$ and the associated Hermitian theta series is a Hermitian modular form of weight $r$ for the full modular group. An explicit description of the isometry classes of these lattices has so far only been obtained for the Gaussian number field $K=\Q(i)$ and $r=4,8$ by Schiemann \cite{S}, where $1$ resp. $3$ isometry classes exist, and for $r=12$ by Kitazume and Munemasa \cite{KM}, where $28$ isometry classes exist. It should be noted that similar results over the Hurwitz order can be found in \cite{Q} and \cite{BN}.

In this paper we derive analogous results for $\Ocal_K$-modules whenever $K=\Q(i\sqrt{3})$ is Eisenstein's number field. There is exactly one isometry class for $r=4$ and $r=8$ and surprisingly there exist just $5$ isometry classes if $r=12$. In the proof we make use of the Niemeier classification and verify which of these lattices have got the structure of an $\Ocal_K$-module. The main idea is that the orthogonal automorphism group of such a lattice must contain an element with minimal polynomial $p(X)=X^2-X+1$, which corresponds to $\frac{1}{2}(1+i\sqrt{3})I$. Moreover we investigate the filtration of cusp forms spanned by the associated Hermitian theta series.

\section{The mass of the Eisenstein lattices}

Throughout the paper let $K=\Q(i\sqrt{3})$ be Eisenstein's number field of discriminant $-3$ with attached Dirichlet character
\[
 \chi_K(n) = \left(\frac{n}{3}\right),\;\; n\in\Z,
\]
and ring of integers
\[
 \Ocal:=\Z+\Z\omega, \quad \omega =(1+i\sqrt{3})/2,
\]
which is Euclidean with respect to the norm $N(a)=a\overline{a}$. We consider $K^r$, $r\in \N$, with the standard Hermitian scalar product 
\[
 K^r\times K^r \to K,\quad (x,y)\mapsto \langle x,y\rangle:=\overline{x}^{tr}y.
\]
Let $\Ucal(r)$ denote the unitary group of $\rank r$.
\begin{definition}
 $\Lambda\subset K^r$ is called an \emph{Eisenstein lattice} of $\rank r$ if there exist linearly independent $b_1,\ldots,b_r\in K^r$ such that 
\begin{enumerate}[(i)]
 \item $\Lambda=\Ocal b_1+\ldots +\Ocal b_r$,
 \item $\det(\langle b_j,b_k\rangle)= (2/i\sqrt{3})^r$,
 \item $\langle \lambda,\lambda\rangle \in 2\Z$ for all $\lambda\in\Lambda$.
\end{enumerate}
Two such lattices $\Lambda,\Lambda'$ are called \emph{isometric} if there exists an \emph{isometry} $U\in\Ucal(r)$ such that
\[
 \Lambda' = U\Lambda.
\]
The (unitary) \emph{automorphism group} of $\Lambda$ is defined by
\[
 \aut (\Lambda):=\{U\in \Ucal(r);\;U\Lambda = \Lambda\}.
\]
The Eisenstein lattices form a genus and its \emph{mass} is
\[
 \mu_r:=\sum_{\Lambda}\frac{1}{\sharp \aut (\Lambda)},
\]
where we sum over representatives of the isometry classes of Eisenstein lattices of $\rank r$. 
\vskip2ex
Considering the symmetric bilinear form
\[
 K^r\times K^r\to \Q,\quad (x,y)\mapsto \re\langle x,y\rangle,
\]
we may also consider $\Lambda$ as an even unimodular $\Z$-lattice of $\rank 2r$ 
\[
 \Lambda_{\Z}= \Z b_1+\Z\omega b_1+\ldots +\Z b_r+\Z\omega b_r.
\]
\end{definition}
This yields an embedding of $\aut \Lambda$ into the (orthogonal) automorphism group $\aut_{\Z}(\Lambda)$ of $\Lambda_{\Z}$. Hence it is clear that Eisenstein lattices exist only if $r$ is a multiple of $4$.

\begin{theorem}\label{theorem 1}
The mass of the genus of Eisenstein lattices of rank $r$, $r\in\N$, $4\mid r$, is given by
\[
 \mu_r = \frac{1}{2^{r-1}\cdot r!} \cdot \prod^{r/2}_{j=1} |B_{2j}\cdot B_{2j-1,\chi_K}|,
\]
where $B_j$ (resp. $B_{j,\chi_K}$) is the (generalized) Bernoulli number. In particular one has
\[
\mu_4 = \frac{1}{2^7\cdot 3^5\cdot 5}, \quad \mu_8= \frac{1}{2^{15}\cdot 3^{10}\cdot 5^2},\quad \mu_{12}=\frac{691\cdot 809\cdot 1847}{2^{22}\cdot 3^{17}\cdot 5^3\cdot 7\cdot 11\cdot 13}\,.
\]
\end{theorem}
\begin{eproof}
Hashimoto and Koseki \cite{HaK} computed the mass $\mu^*_r$ of odd unimodular $\Ocal$-modules of $\rank r$, $4\mid r$, to be
\[
\mu^*_r = \frac{3^{r/2}+1}{2^r\cdot r!}\cdot\prod^{r/2}_{j=1}|B_{2j}\cdot B_{2j-1,\chi_K}|,
\]
where $B_j$ (resp. $B_{j,\chi_K}$) are the (generalized) Bernoulli numbers (cf. \cite{Z}, resp. \cite{F}). Using the same counting argument as in \cite{BN}, Proposition 2.4, we obtain
\[
 \mu_r=\mu^*_r\cdot \frac{c_1}{c_2},
\]
where $c_1$ is the number of maximal isotropic subspaces of the orthogonal $\F_3$-vector space $\F^r_3$ and where $c_2$ is the number of maximal isotropic subspaces of the symplectic $\F_3$-vector space $\F^r_3$. These numbers can be found in \cite{T}, p. 78 resp. p. 174, and are equal to 
\[
 c_1  = \prod^{r/2}_{j=1} (3^{j-1}+1), \quad  c_2  = \prod^{r/2}_{j=1} (3^j+1).
\]
This yields the first claim. The special values are obtained from an explicit calculation of the Bernoulli numbers exploiting the formulas in \cite{Z}.
\qed
\end{eproof}

\section{Classification of Eisenstein lattices}

In order to obtain a classification of the Eisenstein lattices we make use of Theorem 
 \ref{theorem 1}. First of all we use the fact from \cite{DK} that there exists an Eisenstein lattice $\Lambda_4$ of $\rank 4$ given by its Gram matrix
\[
 \begin{pmatrix}
  2I & C   \\
  -C & 2I
 \end{pmatrix}, \quad I=\begin{pmatrix}
			1 & 0  \\
			0 & 1
			\end{pmatrix}, \quad C=\frac{2i}{\sqrt{3}}\begin{pmatrix}
								   1 & 1  \\
								   1 & -1
								   \end{pmatrix}.
\]
An easy computer calculation yields
\[
 \sharp\,\aut (\Lambda_4) = 2^7\cdot 3^5\cdot 5 = \frac{1}{\mu_4}.
\]
In view of 
\[
 \sharp\,\aut(\Lambda_4\oplus \Lambda_4) = 2\cdot(\sharp\,\aut(\Lambda_4))^2 = \frac{1}{\mu_8}
\]
Theorem \ref{theorem 1} implies
\begin{corollary}\label{corollary 1}
There is exactly one isometry class of Eisenstein lattices of $\rank 4$ and of $\rank 8$. Representatives are $\Lambda_4$ and $\Lambda_4\oplus \Lambda_4$.
\end{corollary}
In the case of $\rank 12$ we apply the Niemeier classification of even unimodular $24$-dimensional $\Z$-lattices (cf. \cite{CS}). We use the fact that $\omega I$ is always contained in the automorphism group of an Eisenstein lattice. Its minimal polynomial over $\Q$ is $p(X)= X^2-X+1$.

\begin{proposition}\label{proposition 1}
Each matrix $B\in\Z^{2r\times 2r}$ with minimal polynomial $p(X) = X^2-X+1$ is in ${\rm GL}(2r;\Z)$ conjugate to
\[
 \diag(A,\ldots,A),\quad A=\begin{pmatrix}
                            0 & -1  \\
			    1 & 1
                           \end{pmatrix}.
\]
\end{proposition}
\begin{eproof}
Following Newman \cite{N}, p. 54 and Theorem III.12, the claim basically follows from the fact that $\Ocal = \Z[\omega]$ has class number $1$.
\qed
\end{eproof}
\vskip2ex

Now we are going to investigate which Niemeier lattices possess an automorphism with minimal polynomial $p(X)=X^2-X+1$.

\begin{lemma}\label{lemma 1}
Let $L$ be a $\Z$-lattice, which has a decomposition into orthogonally indecomposable  non-isometric lattices
\[
 L = R^{n_1}_1 \oplus \ldots \oplus R^{n_s}_s.
\]
Then one has  \\[1ex]
{\rm a)} $\aut(L) = \text{\Large $\times$}^s_{j=1}\aut(R_j^{n_j})$.  \\[1ex]
{\rm b)} $\aut(R^{n_j}_j) = \{\diag(\phi_1,\ldots,\phi_{n_s})\sigma_j;\;\phi_k\in\aut(R_j),\sigma_j\in S_{n_j}\}$.
    \\[1ex]
{\rm c)} If $\phi\in\aut(L)$ exists with minimal polynomial $p(X) = X^2-X+1$, then  $\phi$ belongs to the kernel of the group epimorphism $\aut(L) \to S_{n_1}\times \ldots\times S_{n_s}$ in {\rm a)} and {\rm b)}.
\end{lemma}
\begin{eproof}
a), b) The results are well-known (cf. \cite{Kn}).   \\
c) Let $P_{\sigma}$ be the integral orthogonal $n\times n$ matrix associated to $\sigma\in S_n$ in the $j$-component. Then $p(\phi)=0$ shows that $P_{\sigma}^2-P_{\sigma}$ is diagonal. Hence $\sigma=id$ follows.
\qed
\end{eproof}
\vskip2ex
Thus we have reduced the problem to orthogonally indecomposable $\Z$-lattices. In the case of root lattices (cf. \cite{CS}, Chapter 4) we obtain the well-known 

\begin{lemma}\label{lemma 2}
a) If $n\geq 2$ one has
\[
 \aut(A_n)\cong S_{n+1}\times \{\pm 1\}.
\]
b) If $n\neq 4$ one has
\[
 \aut(D_n)\cong S_n \ltimes \{\pm 1\}^n.
\]
\end{lemma}
Considering our situation we obtain

\begin{corollary}\label{corollary 2}
If
\[
 \Lambda = A_n,\;\;n\neq 2 \quad \text{or}\quad \Lambda = D_n,\;\;n\neq 4 \quad \text{or}\quad \Lambda=E_7
\]
then $\aut(\Lambda)$ does not contain an automorphism with minimal polynomial $p(X) = X^2-X+1$.
\end{corollary}
\begin{eproof}
The claim is clear for $A_1=\Z$. Suppose that $\phi\in\aut(\Lambda)$ with minimal polynomial $p(X)=X^2-X+1$ exists. Hence $\phi^3=-id$ follows. Considering $A_n$, $n\geq 3$, and the action of $S_{n+1}$ on $A_n$ we conclude that $A_n$ contains a sublattice of $\rank 3$ which is $\phi$-invariant. As $3$ is odd $\phi$ must possess a real eigenvalue which contradicts $p(X)=X^2-X+1$.   \\
$\aut(D_n)$, $n=1,2$, does not contain an element of order $3$. Hence let $n=3$ or $n\geq 5$ and $\Lambda=D_n$. Again we can conclude that $D_n$ contains a sublattice of $\rank 3$ which ist $\phi$-invariant. Hence the same argument as above yields the contradiction. \\
In view of $\rank E_7=7$ we may proceed in the same way as in the last case.
\qed
\end{eproof}
\vskip2ex
A verification yields
\begin{lemma}\label{lemma 3}
The lattices $A_2,\, D_4,\, E_6,\, E_8$ and the Leech lattice $L_0$ possess an automorphism with minimal polynomial $p(X)=X^2-X+1$. This automorphism is unique up to conjugacy.
\end{lemma}

\begin{eproof}
The existence and the uniqueness are verified by a direct computation for $A_2$ and $D_4$. Considering $E_6$, $E_8$ and $L_0$ the result follows from \cite{CC}.
\qed
\end{eproof}
\vskip2ex
Hence we get

\begin{theorem}\label{theorem 2}
There are exactly $5$ isometry classes of Eisenstein lattices of $\rank 12$. The underlying $\Z$-lattices, the orders of the unitary automorphism groups and their index in the orthogonal automorphism groups are given by the following table
\begin{center}
\begin{tabular}{>{$}c<{$}|>{$}c<{$}|>{$}c<{$}}
 \text{root system of }\Lambda & \sharp\,\aut (\Lambda) & {\rm [}\aut_{\Z}(\Lambda):\aut(\Lambda){\rm ]}\\  
 \hline
 3E_8 & 2^{22}\cdot 3^{16}\cdot 5^3 & 2^{21}\cdot 5^3\cdot 7^3 \\
 4E_6 & 2^{16}\cdot 3^{17} & 2^{16}\cdot 5^4  \\
 6D_4 & 2^{21}\cdot 3^9\cdot 5 & 2^{19} \\
 12A_2 & 2^7 \cdot 3^{15}\cdot 5\cdot 11 & 2^{12} \\
 \emptyset & 2^{14}\cdot 3^8\cdot 5^2\cdot 11\cdot 13 & 2^8\cdot 3\cdot 5^2\cdot 7\cdot 23 
\end{tabular}
\end{center}
\vspace{1ex} 
\end{theorem}

\begin{eproof}
Apply Lemma \ref{lemma 3} and the Niemeier classification in \cite{CS}. After conjugation we may assume that 
\[
B=\diag(A,\ldots,A), \quad A=\begin{pmatrix}
				0 & -1 \\ 1 & 1 
				\end{pmatrix},
\] 
belongs to $\aut(\Lambda)$ due to Proposition \ref{proposition 1}. Using
\[
 A\Z^2\cong \Ocal
\]
we have got an $\Ocal$-structure on $\Lambda$. Moreover we have 
\[
 \aut(\Lambda) \cong \{U\in\aut_{\Z}(\Lambda);\;UB=BU\}.
\]
Now $\sharp \aut(\Lambda)$ and $[\aut_{\Z}(\Lambda):\aut(\Lambda)]$ are computed explicitly using \cite{PS}. Using Theorem \ref{theorem 1} we conclude that we have already found all the isometry classes.
\qed
\end{eproof}

\begin{remark}\label{remark 1}
a) Explicit examples of the lattices in Theorem 2 given by their Gram matrices can be found in \cite{H1}.\\[1ex]
b) Theorem \ref{theorem 1} yields
\[
 \mu_{16} = \frac{13\cdot 47\cdot 419\cdot691\cdot 809\cdot 1847\cdot 3617\cdot 16519}{2^{31}\cdot 3^{22}\cdot 5^4\cdot 11\cdot 17}\;\approx\;0,002.
\]
So far we have constructed $5$ isometry classes of Eisenstein lattices $\Lambda_j'$ of $\rank 16$, namely the orthogonal sum from those in Theorem \ref{theorem 2} with $\Lambda_4$. But they have got only a small proportion of the whole mass 
\[
 \sum^5_{j=1} \frac{1}{\sharp \aut(\Lambda_j')}\bigg/ \mu_{16} < 5\cdot 10^{-14}.
\]
c) The Eisenstein lattices associated with $3E_8,\,6D_4$ and $L_0$ have also got the structure of a module over the Hurwitz order (cf. \cite{BN}). The quaternionic automorphism group has the index $3^{12}$ resp. $3^6$ resp. $2\cdot 3^5\cdot 11$ in the unitary automorphism group $\aut(\Lambda)$. 
\vspace{1ex}\\
d) Analogous results for other imaginary quadratic number fields and $r=4,8$ are achieved in \cite{H2}.
\end{remark}

\section{The Hermitian theta series}

We consider the Hermitian half-space of degree $n$
\[
 \Hcal_n = \{Z\in\C^{n\times n};\;\frac{1}{2i}(Z-\overline{Z}^{tr}) > 0\}.
\]
The Hermitian modular group
\[
 \Gamma_n:=\{M\in {\rm SL}(2n;\Ocal);\;MJ\overline{M}^{tr}=J\},
\quad J= \begin{pmatrix}
          0 & I   \\
	  -I & 0
         \end{pmatrix},
\]
acts on $\Hcal_n$ by the usual fractional linear transformation. The space $[\Gamma_n,r]$ of Hermitian modular forms of degree $n$ and weight $r$ consists of all holomorphic $f:\Hcal_n \to \C$ satisfying
\[
 f(M\langle Z\rangle) = \det(CZ+D)^r \cdot f(Z)\quad \text{for all}\quad M=\begin{pmatrix}
	A & B  \\  C & D
        \end{pmatrix}\in\Gamma_n
\]
with the additional condition of boundedness for $n=1$ (cf. \cite{B}). The subspace $[\Gamma_n,r]_0$ of cusp forms is defined by the kernel of the Siegel $\phi$-operator
\[
 f\mid \phi \equiv 0.
\]
\begin{proposition}\label{proposition 2}
{\rm (\cite{CR})} Let $\Lambda$ be an Eisenstein lattice of $\rank r$ with Gram matrix $H$. Then the associated Hermitian theta series
\begin{align*}
 \Theta^{(n)} (Z,H): & =\sum_{G\in \Ocal^{r\times n}} e^{\pi i\trace(\overline{G}^{tr}HG\cdot Z)}, \;\; Z\in \Hcal_n,  \\
& = \sum_{T\geq 0} \sharp(H,T) e^{\pi i \trace (TZ)},
\end{align*}
with the Fourier coefficients
\[
\sharp(H,T):= \sharp\{G\in \Ocal^{r\times n};\;\overline{G}^{tr}HG = T\},
\]
belongs to 
\[
 [\Gamma_n,r].
\]
\end{proposition}

Just as in the case of Siegel modular forms we obtain the analytic version of Siegel's main theorem involving the Siegel Eisenstein series in $[\Gamma_n,r]$ (cf. \cite{B}, \cite{K}).

\begin{corollary}\label{corollary 3}
Let $H_1,\ldots,H_s$ be Gram matrices of the lattices $\Lambda_1,\ldots,\Lambda_s$, which are representatives of the isometry classes of  Eisenstein lattices of $\rank r$. If $r>2n$ one has 
\[
 \frac{1}{\mu_r}\sum^s_{j=1}\frac{1}{\sharp\aut(\Lambda_j)}\Theta^{(n)}(Z,H_j) = E^{(n)}_r(Z) = \sum_{\left(\begin{smallmatrix}
A & B \\ C & D  \end{smallmatrix}\right) : \left(\begin{smallmatrix}
\ast & \ast \\ 0 & \ast \end{smallmatrix}\right)\big\backslash \Gamma_n} \det(CZ+D)^{-r}.
\]
\end{corollary}
If $r=4$ or $r=8$ the Eisenstein series therefore coincides with a Hermitian theta series. Let $H_1,\ldots,H_5$ be fixed Gram matrices of the Eisenstein lattices associated with $3E_8,4E_6,6D_4,12A_2$ resp. $L_0$. Then we can describe the filtration analogous to \cite{NV}.  

\begin{theorem}\label{theorem 3}
a) Let 
\begin{align*}
F^{(n)}(Z):=\Theta^{(n)} & (Z,H_1)-30\Theta^{(n)}(Z,H_2)+135 \Theta^{(n)}(Z,H_3) -160\Theta^{(n)}(Z,H_4)    \\
 & +54 \Theta^{(n)}(Z,H_5)\in [\Gamma_n,12].
\end{align*}
Then $F^{(4)}\not\equiv 0$ and either $F^{(3)}\equiv 0$ or $F^{(3)}$ is a non-trivial cusp form.  \\
b) The form
\[
 G^{(3)}(Z) = \Theta^{(3)}(Z,H_2) -6\Theta^{(3)}(Z,H_3)+8\Theta^{(3)}(Z,H_4)-3\Theta^{(n)}(Z,H_5)
\]
is a non-trivial cusp form of degree $3$.  \\
c) The form
\[
 H^{(2)}(Z) = \Theta^{(2)}(Z,H_3)-2\Theta^{(2)}(Z,H_4)+\Theta^{(2)}(Z,H_5)
\]
is a non-trivial cusp form of degree $2$.  \\
d) The form
\[
 J^{(1)}(Z) = \Theta^{(1)}(Z,H_4) -\Theta^{(1)}(Z,H_5)
\]
is a non-trivial cusp form of degree $1$.
\end{theorem}

\begin{eproof}
The functions above are modular forms of weight $12$ due to Proposition \ref{proposition 2}.  \\
a) The modular form $F^{(4)}$ does not vanish identically, since this is already true for the restriction to the Siegel half-space, where $F^{(4)}$ is the analogous linear combination of the Siegel theta series. In this case we may use the table of Fourier coefficients in \cite{BFW}.  \\
We know from \cite{DK} that 
\[
 [\Gamma_2,12] = \C E_4^{(2)^3}+\C E_6^{(2)^2}+\C E_{12}^{(2)}.
\]
Calculating a few Fourier coefficients we obtain $F^{(2)}\equiv 0$.   \\
b) The analogous computations show that $G^{(2)}\equiv 0$. If $n=3$ the Fourier coefficient of 
\[
 \begin{pmatrix}
  4 & 4\omega/3 & 1+\omega/3  \\
  4\overline{\omega}/3 & 4 & -2\omega/3  \\
  1+\overline{\omega}/3 & -2\overline{\omega}/3 & 4
 \end{pmatrix}
\]
is $5184$, hence $G^{(3)}\not\equiv 0$.  \\
c)We have
\[
 H^{(2)} = c\cdot \left(E_{12}^{(2)}-\frac{441}{691}\left(E_4^{(2)}\right)^3-\frac{250}{691}\left(E_6^{(2)}\right)^2\right) \in[\Gamma_{2},12]_0,\;\;c\neq 0.
\]
d) This claim is clear as we deal with elliptic modular forms and elliptic theta series, hence
\[
 J^{(1)} = 720 \Delta,
\]
where $\Delta$ is the normalized discriminant.
\qed
\end{eproof}

\begin{remark}\label{remark 2}
a) Just as in the case of quaternionic modular forms (cf. \cite{HK}) there is no analog of the Schottky form over $K$.   \\
b) A computation of thousands of Fourier coefficients leads to the conjecture that $F^{(3)}\equiv 0$. In this case the space of cusp forms of weight $12$ spanned by Hermitian theta series is one dimensional in the cases of degree $n=0,1,2,3,4$. One knows that 
\[
 [\Gamma_4,12]_0\neq \{0\}
\]
from \cite{I}. Hence it is conjectured that $F^{(4)}$ is a Hermitian Ikeda lift.\\
c) All the theta series associated with Eisenstein lattices of $\rank\leq 12$ are symmetric Hermitian modular forms, i.e. $f(Z^{tr})=f(Z)$, because the lattice $\overline{\Lambda}$ is isometric to $\Lambda$.
\end{remark}

\vspace{2ex}

\renewcommand\refname{References}

\end{document}